\documentclass[11pt,leqno]{amsart}

\usepackage[english]{babel}
\usepackage{a4wide}
\usepackage[utf8]{inputenc}
\usepackage{amsmath}
\usepackage{amssymb}
\usepackage{amsthm}
\usepackage{amsfonts}
\usepackage{amsmath}
\usepackage{latexsym}
\usepackage{ifthen, enumitem, mathrsfs, euscript,anyfontsize}
\usepackage{graphicx}
\usepackage[table]{xcolor}
\usepackage{listings}
\usepackage{url}
\usepackage{amsfonts}
\usepackage{enumerate}
\usepackage{blindtext}
\usepackage{hyperref}
\usepackage{xcolor}
\usepackage{ bbold } 
\usepackage{color}

\newcommand{\one}{\mathbb 1}

\newtheorem*{theorem}{Theorem}

\newtheorem*{lemma}{Lemma}

 \everymath{\displaystyle}
\def\Z{\mathbb{Z}}

\def\R{\mathbb{R}}

\def\P{{\mathbb P}}

\def\cP{{\EuScript P}}
\def\cQ{{\EuScript Q}}
\def\cT{{\EuScript T}}
\renewcommand{\le}{\leqslant }

\renewcommand{\leq}{\leqslant }
\renewcommand{\geq}{\geqslant }
\renewcommand{\epsilon}{\varepsilon}
\newcommand{\mmod}[1]{\,({\rm mod\,}#1)}

\newcommand{\e}{{\rm e}}
\newcommand{\sset}{\smallsetminus}
\makeatletter
\def\@settitle{\begin{center}%
  \baselineskip14\p@\relax
    \normalfont\LARGE
  \@title
  \end{center}%
}
\makeatother
\begin{document}
\let\MakeUppercase\relax 
\title{On simply normal numbers \goodbreak with digit dependencies}
\author{\large Ver\'onica Becher, Agust\'in Marchionna,  G\'erald Tenenbaum}
\date{April 12, 2023}

\subjclass[2020]{Primary 11K16, 11N60; Secondary 11N56}
\keywords{normal numbers, Toeplitz sequences, discrepancy, additive and multiplicative functions}
\maketitle

\begin{abstract}
Given an integer $b\geq 2$ and a set $\cP$ of prime numbers, 
the set   $\cT_\cP $ of  Toeplitz  numbers comprises 
all elements of $[0,b[$  whose digits $(a_n)_{n\geqslant 1}$ in the base-$b$ expansion 
satisfy $a_n=a_{pn}$ for all  $p\in \cP$ and $n\geq 1$. 
Using a completely additive arithmetical function, we construct a number in~$\cT_\cP$ 
that is simply Borel normal if, and only if, 
 $\textstyle \sum_{p\not \in\cP} 1/p=\infty$.
We then provide an effective  bound for the discrepancy.
\end{abstract}
\medskip

\bigskip

Let $\mathbb P$ denote the set of prime numbers, and  let $\cP\subset {\mathbb P}$.
Following Jacobs and Keane's  definition of Toeplitz sequences  in~\cite{JacobsKeane69},
we define the set  $\cT_{\mathcal P}$  of   \emph{Toeplitz numbers} 
  as  the set of all real numbers $\xi\in[0,b[$ whose base-$b$ expansion 
$\textstyle \xi=\sum_{n\geq 1} a_n/b^{n}$ 
satisfies  
\[
a_n = a_{np}\qquad (n\geqslant 1,\,p\in\cP).
\]
For example, $0.a_1 a_2 a_3 \dots$ 
 is a Toeplitz number for $\cP=\{2,3\}$ if, for every   $n \geq 1$, we have
\[
a_{n} = a_{2n} = a_{3n}.
\]
Then, $ a_1, a_5, a_7, a_{11}, \dots$ are independent
 while $a_2,a_3,a_4,a_6,\dots$ are 
completely determined by earlier digits.

As defined by \'Emile Borel, a real number is called {\em simply normal} to the integer base
$b\geq 2$ if every possible digit in $\Z/b\Z$
 occurs in its $b$-ary
expansion  with the same asymptotic frequency~$1/b$.  A real number is said to be {\em normal} to the base $b$ if it is simply
normal to all the bases~$b^j$, $j\geqslant 1$.  
Borel  proved that, with respect to the Lebesgue measure, almost all numbers are  normal to all integer bases 
at least equal to $2$. 
For a presentation of the theory of normal numbers see for
example~\cite{Bugeaud2012,KN1974}.

In~\cite{ABC}, Aistleitner, Becher and Carton considered the notion of 
 Borel normality under the assumption of 
dependencies between the digits of the expansion.
Thus~\cite[th.~1]{ABC} 
states that, given any integer base~$b \geq 2$  
and any finite subset~$\cP$ of the primes,
almost all  numbers, with respect to the 
 uniform probability measure on $\cT_\cP$,
are normal to the base~$b$.
In  the particular case  $\cP = \{2\}$, they 
 show  \cite[th.~2]{ABC} 
that almost all   elements in $\cT_\cP$ (still with respect to the uniform  measure on $\cT_\cP$)
are  normal to all integer bases greater than or equal to~$2$.
For $\cP=\{2\}$, a construction of an explicit  number in $\cT_\cP$ that is normal to the base $2$
 appears in~\cite{BCH}.  This construction can be generalized to any integer base $b$ and any singleton~$\cP$.  

Let $\Omega_{\cP}$ denote the completely additive arithmetical function defined by 
$\Omega_{\cP}(p)=\one_{(\P\smallsetminus\cP)} (p)$.
Then,  $\Omega_\cP(n)$ is the sum of the exponents in the canonical factorization of~$n$ of those prime factors  that do {\em not} belong to  $\cP$. For $n\geqslant 1$ and $b\geqslant 2$, let $a_n=a_{n,b}$ denote the representative of the congruence class $\Omega_\cP(n)+b\Z$ lying in $[0,b[$.
Thus, given $b\geq 2$, the real number
\begin{equation}
\label{defxiP}
 \xi_\cP=\sum_{n\geqslant  1} a_n/b^{n}
 \end{equation}
is a   an element of  $\cT_\cP$.

Motivated by the question posed in~\cite{ABC} on  how to exhibit 
a normal number in $\cT_\cP$ for any set  $\cP$ of primes,  
we construct in this note  simply normal numbers  for arbitrary bases and a large family of sets $\cP$.
\begin{theorem}\label{thm:main}
Let $\cP\subset \P$, $\cQ:=\P\sset\cP$,  and let $b$ be an integer $\geq 2$.
The number $\xi_\cP$
is  simply normal to the base $b$ if, and only if,  
\begin{equation}
\label{cns}\sum_{p\in\cQ} 1/p=\infty.
\end{equation}
Moreover, defining, for $0\leqslant k<b$,
\[
\epsilon_{N,k}:=
\left|\frac1N |\{1\le n \le N : a_n=k\}| -\frac{1}{b}
\right|,\quad E(N):=\sum_{p\leq N,\,p\in\cQ} \frac1p \quad(N\geqslant 1),
\]
we have
\begin{equation}
\label{eff}
\epsilon_{N,k} \ll 
\e^{-2E(N)/9b^2}.
\end{equation}
\end{theorem}

Our proof rests on the following auxiliary result where we use the traditional notation $\e(u):=\e^{2\pi i u}$ $(u\in\R)$.
\begin{lemma}\label{average}
 Let $\cP\subset \P$  and let $b$ be an integer $\geq 2$.
The number  $\xi_\cP$
is simply normal to the base~$b$ if, and only if,
\begin{equation}
\label{crit}\frac{1}{N} \sum_{n\leq N}  \e\big({a\Omega_{\cP}(n)/b}\big) =o(1) \quad(a=1,2,\ldots b-1,\,N\to\infty).
\end{equation}
\end{lemma}

\begin{proof}
The necessity of the criterion is clear.
We show the sufficiency.
Define 
\[
b_{k,N}=\frac1N\big|\{1\leq n \leq N : a_n=k\}\big|\quad(0\leqslant k<b,\,N\geqslant 1).
\]
Then
\begin{equation}
\label{evalb}
b_{k,N}=\frac1{bN}\sum_{0\leqslant a<b}\e(-ak/b)\sum_{1\leqslant n\leqslant N}\e\big(a\Omega_\cP(n)/b\big)=\frac1b+o(1)
\end{equation}
since by \eqref{crit} all inner sums with $a\neq0$ contribute $o(N)$. 
\end{proof}
We may now embark on the proof of the Theorem.
Let 
\begin{equation*}
S(N; a/b) := \sum_{n \leq N} \e\big(a\Omega_\cP(n)/b\big)\quad(a\in\Z,\,b\geqslant 2,\,N\geqslant 1).
\end{equation*}
 We aim at necessary and sufficient conditions that ensure $S(N, a/b) = o(N)$ as $N \rightarrow +\infty$, and seek effective upper bounds for $S(N;a/b)$ when such conditions are met.
 
Whenever  $a$ and $b$ are coprime, $b \geq 2$ and $|a| \leq b/2$, we may apply
\cite[cor.\thinspace2.4(i)]{T2017} with   $r=1$, $z=\e(a/b)$, $\vartheta=2\pi a/b$ and $\kappa = 1$. Using \cite[(7.4)]{T2017}, from which the bound \cite[(2.19)]{T2017} is actually derived, this yields 
\begin{equation*} \label{effective}
S(N; a/b) \ll Ne^{-2a^2E(N)/(9b^2)}.
\end{equation*}
So, if \eqref{cns} holds, then the above  lemma implies that
$\xi_\cP$ is simply normal to the base $b$.
Notice that $\{a\in\Z:|a|\leqslant \tfrac12b\}$ describes a complete set of residues $\mmod b$. 
The effective bound \eqref{eff} is then provided by \eqref{evalb}.

If, on the contrary,  condition \eqref{crit} fails, we  apply \cite[cor.\thinspace2.2]{T2017}, which is an effective version of a result of Delange---see \cite[th. III.4.4]{Tbook}. We have 
\begin{equation}
\label{majconv}    \sum_{p \in \cQ,~p \leq N} \frac{\log p}{p} \ll \eta_N\log N
\end{equation}
for some $\eta_N \to 0$.
A possible choice is 
\[
\eta_N:=\min_{1\leq z\leq N} \bigg(\frac{\log z}{\log N}+\sum_{p\in\cQ,\,p>z} \frac1{p}\bigg).
\]
The validity of \eqref{majconv} is then obtained by bounding $\log p$ by $\log z$ if $p\leqslant z$ and by $\log N$ otherwise. That $\eta_N=o(1)$ follows by noticing that the last sum tends to 0 as $z\to\infty$.
Then we get 
\begin{equation*}
    S(N;a/b)=\frac{N}{\log N} 
    \left( 
    \prod_p \sum_{p^\nu \le N}
    \frac{e(\nu a\Omega_{\cP}(p)/b)}{p^\nu}
+O\left(\eta_N^{1/8}\e^{E(N)}+\frac{\e^{E(N)}}{(\log N)^{1/12}}\right)  
    \right),
\end{equation*}
where we are picking the corresponding values from \cite[cor.~2.2]{T2017} as $a=1/8$, $b=1/12$, and~$\varrho=1$.\par 
 To prove that 
\begin{equation}
\label{minS}
S(N,a/b)\gg N,
\end{equation}
it hence suffices to show  that
\[
\log N\ll\prod_p \sum_{p^\nu \le N} 
    \frac{e(\nu a\Omega_{\mathcal P}(p)/b)}{p^\nu}=\prod_{p\in\cQ}\frac{1-\e(\nu_pa/b)/p^{\nu_p}}{1-\e(a/b)/p}\prod_{p\in\cP}\frac{1-1/p^{\nu_p}}{1-1/p}
,
\]
where we have put $\nu_p:=1+\lfloor{(\log N)/\log p}\rfloor$, so that $p^{\nu_p}\geqslant N$. Now the double product above is clearly
$$\sim \sigma_N:=\prod_{p\leqslant N}\frac{1}{1-1/p}\prod_{p\in\cQ}\frac{1-1/p}{1-\e(a/b)/p}\cdot$$
Since the general factor of the last product equals $1+\{\e(a/b)-1\}/p+O(1/p^2)$, we deduce from the convergence of $\textstyle\sum_{p\in\cQ}1/p$ and Mertens' formula that $\sigma_N\sim c\log N$ for some $c\neq0$. 
This yields \eqref{minS} as required.

\bibliographystyle{plain}
\bibliography{simply}

\footnotesize
\bigskip
\noindent
Ver\'onica Becher\\
Departamento de Computaci\'on, Facultad de Ciencias Exactas y Naturales, Universidad de Buenos Aires   e ICC CONICET 
\\
Pabell\'on 0, Ciudad Universitaria, C1428EGA Buenos Aires, Argentina
\\
{\tt vbecher@dc.uba.ar}
\medskip
\medskip

\noindent
Agust\'in Marchionna
\\
Departamento de Computaci\'on, Facultad de Ciencias Exactas y Naturales, Universidad de Buenos Aires 
\\Pabell\'on 0, Ciudad Universitaria, C1428EGA Buenos Aires, Argentina
\\
{\tt agusmarchionna1998@gmail.com}
\medskip
\medskip

\noindent
G\'erald Tenenbaum
\\
Institut \'Elie Cartan, Universit\'e de Lorraine
\\
BP 70239
\\
54506 Vand\oe uvre-l\`es-Nancy Cedex France
\\
{\tt gerald.tenenbaum@univ-lorraine.fr}
\medskip

\end{document}